# Comparative study of Three Numerical Schemes for Fractional Integro –Differential Equations

Kamlesh Kumar*, Rajesh K. Pandey[1]* and Shiva Sharma*

*Department of Mathematical Sciences, Indian Institute of Technology (BHU) Varanasi, UP India

## Abstract:

This paper presents a comparative study three numerical schemes such as Linear, Quadratic and Quadratic-Linear scheme for the fractional integro-differential equations defined in terms of the Caputo fractional derivatives. The error estimates of the respective approximations are also established. Numerical tests of the discussed schemes show that all schemes work well, and when the number of terms approximating the solution are increased, the desired solution is achieved. The accuracy of the numerical schemes with respect to the step size $h$ is analyzed and illustrated through various tables. Finally, comparative performances of the schemes are discussed.

## 1. Introduction

Fractional calculus is a branch of mathematics that deals with generalization of the well-known operations of differentiations and integrations to arbitrary orders. Many papers on fractional calculus have been published during the second half of the 20th century. Recent developments of fractional calculus involve the real world applications in science and engineering such as viscoelasticity [1-3], bioengineering [4], biology [5] and more can be found in [6-8]. Two main differences between fractional calculus and classical calculus are, 1) fractional integrals and derivatives are nonlocal, and 2) in the limit results obtained from fractional calculus coincide with those obtained from classical calculus.  This makes fractional calculus to be richer than classical calculus. Many real world problems are being modelled using fractional derivative and integral terms and such equations are known as the fractional integro-differential equations. Fractional integro-differential equations arise in the areas of signal processing [8],  mechanics [9], econometrics [10], fluid dynamics [7], nuclear reactor dynamics, acoustic waves [11] and electromagnetics [12] etc. There have been much developments on the analytical and numerical methods for solving the fractional integro-differential equations in recent years. Some of them are described as follows: In [13] Saadatmandi and Dehghan applied the Legendre collocation method for solving the fractional integro-differential equations. Rawashdeh [14] presented the spline collocation scheme for solving fractional integro-differential equations. In [15-16], the application of the collocation method has been extended for solving linear and nonlinear fractional integro-differential equations. Galerkin method, and wavelet Galerkin method based numerical methods are studied by the authors respectively in [17] and [18] for fractional integro-differential equations. In [19-20], authors presented the second kind Chebyshev wavelet based approximation and cosine and sine wavelet method respectively for fractional integro-differential equations. Some other methods such as least squares method [21], Tau approximation methods [22-23], Chebyshev pseudo spectral method [24], hybrid collocation method [25], Sinc-Legendre collocation method [26], Legendre spectral element method[27], meshless method [28], and operational matrix method [29] are discussed by the authors for solving fractional integro-differential equations and fractional differential equations in recent years. More recently in 2015, Tohidi and co-authors [30] presented the Euler function based operational matrix approach for solving the fractional integro-differential equations.

In this paper, we present three numerical schemes such as Linear, Quadratic and Quadratic-linear for the fractional integro-differential equations defined in terms of the Caputo derivatives and also discuss the comparative performances of these schemes. These schemes are based on the discretization of the domain into subdomains and then the unknown function is approximated between subdomain using linear and quadratic interpolating polynomials. These approximations schemes are also useful when the analytical integrations of the expression containing the fractional derivative and integral terms are rare or difficult to compute. And thus the numerical schemes become important. The significant contributions in development of these schemes for fractional derivatives are discussed in [31] where the idea of discretized fractional calculus with linear multistep method is presented. Kumar and Agrawal [32] presented the numerical scheme for the fractional differential equations defined in terms of the Caputo derivatives using the quadratic approximations. In 2006, Odibat [33] discussed the algorithms named modified Trapezoidal rules for the approximations of the fractional integral and the Caputo fractional derivatives and also discussed the error estimates. Further, Agrawal [34] discussed the finite element approximation for the fractional variational problem. Later, Pandey and Agrawal [35] presented a comparative study of linear, quadratic and quadratic-linear approximations to the fractional variational problems defined in terms of the generalized derivatives. In [36], the authors presented the higher order approximations to the Caputo derivatives and then applied to solve the advection-diffusion equations. In 2015, Dehghan et al. [37] discussed the high order difference schemes together with Galerkin spectral technique for numerical solution of multi-term time fractional partial differential equations. These recent works motivated the authors to present and discuss the comparative study of Linear, Quadratic and Quadratic-linear schemes for fractional integro-differential equations. Here these schemes are discussed and compared only for linear fractional integro-differential equations. However, similar approach may also be applied for nonlinear case in future. The Quadratic-linear scheme is proposed as combination of the Linear scheme and Quadratic scheme for fractional integro-differential equations. These schemes and their comparative study are analyzed on the illustrative test examples. Further, error estimates of the presented schemes are established. Numerical convergence order of all three schemes are calculated and presented through tables.

## 2. Preliminaries and Statement of Problem

Here, we present some definitions of fractional integral and fractional derivative and define the problem.

*Definition 2.1*: The Riemann–Liouville fractional integration of order $\alpha > 0$ is defined as,

$$(I^\alpha f)(t) = \frac{1}{\Gamma(\alpha)} \int_0^t (t-\tau)^{\alpha-1} f(\tau) d\tau , \qquad (1)$$

$$(I^0 f)(t) = f(t). \qquad (2)$$

And its fractional derivative known as Riemann-Liouville fractional derivative of order $\alpha > 0$ is defined as,

$$(D^\alpha f)(t) = \left( \left( \frac{d}{dt} \right)^n (I^{n-\alpha} f) \right)(t), \qquad n-1 < \alpha \leq n,$$

where $n$ is an integer.

Further, another definition of fractional derivative is introduced and known as Caputo derivative, is defined below.

*Definition 2.2:* The Caputo definition of the fractional differential operator $D^\alpha$ of a function $f(t)$ is defined as,

$$D^\alpha f(t) = \begin{cases} \frac{1}{\Gamma(n-\alpha)} \int_0^x \frac{f^m(\tau)}{(x-\tau)^{\alpha-m+1}} \, d\tau, & m-1 < \alpha < m \\ \frac{d^m}{dx^m} f(t), & \alpha = m \end{cases}, \qquad (3)$$

where, $m$ is an integer. For the detailed study of the fractional calculus, the readers are referred to [6-7].

Now, we consider the general fractional integro-differential equation defined as,

$$D^\alpha \varphi(x) = f(x) + \int_0^x K(x,\tau)\varphi(\tau)d\tau , \ 0 \leq x,\tau \leq 1, \ 0 < \alpha < 1, \qquad (4)$$

with the following supplementary condition $\varphi(0) = \delta$, where $D^\alpha \varphi(x)$ indicates the $\alpha-$th order Caputo derivatives of $\varphi(x)$, and $f(x), K(x,\tau)$ are known functions. The numerical schemes for this problem are presented in the upcoming section.

## 3. Numerical Schemes

In this section, three approximation schemes namely Linear, Quadratic and Quadratic-Linear schemes are presented and discussed for the fractional integro-differential equation defined by Eq. (4). These schemes are denoted as S1, S2 and S3 respectively and are described as follows:

### 3.1 The Linear Scheme (S1)

In the Linear scheme, the domain is divided into $n$ subintervals $[x_j, x_{j+1}]$ and then the unknown function is approximated as the linear interpolation function into each subintervals with the uniform step size $h = \frac{1}{n}$ such that the node points are $x_k = kh$, $k = 0, 1, 2 \ldots n$. And denote the values of functions $\varphi(x), f(x)$ at the node $x_k$ by $\varphi_k = \varphi(x_k)$, and $f_k = f(x_k)$. The fractional integro-differential equation defined by Eq.(4) takes the form for $0 < \alpha < 1$,

$$\frac{1}{\Gamma(1-\alpha)} \int_0^{x_k} (x_k - \tau)^{-\alpha} \varphi'(\tau) d\tau = f(x_k) + \int_0^{x_k} K(x_k, \tau) \varphi(\tau) d\tau. \tag{5}$$

The left side of Eq. (5) can be approximated using the initial condition, $D^\alpha \varphi(x_0) = 0$,

$$\frac{1}{\Gamma(1-\alpha)} \int_0^{x_k} (x_k - \tau)^{-\alpha} \varphi'(\tau) d\tau,$$

$$= \frac{1}{\Gamma(1-\alpha)} \sum_{j=0}^{k-1} \int_{x_j}^{x_{j+1}} (x_k - \tau)^{-\alpha} \varphi'(\tau) d\tau,$$

$$= \frac{1}{\Gamma(1-\alpha)} \sum_{j=0}^{k-1} \int_{x_j}^{x_{j+1}} (x_k - \tau)^{-\alpha} \frac{\varphi_{j+1} - \varphi_j}{h} d\tau,$$

or,

$$\frac{1}{\Gamma(1-\alpha)} \int_0^{x_k} (x_k - \tau)^{-\alpha} \varphi'(\tau) d\tau = \sum_{j=0}^{k} A(k,j) \varphi_j, \tag{6}$$

where,

$$A(k,j) = \frac{h^{-\alpha}}{\Gamma(1-\alpha)} \begin{cases} (k-1)^{1-\alpha} - (k)^{1-\alpha}, & j = 0. \\ (k-j-1)^{1-\alpha} - 2(k-j)^{1-\alpha} + (k-j+1)^{1-\alpha}, & 1 \le j \le k. \\ 1 & j = k \end{cases} \tag{7}$$

Similarly the integral part on right side of Eq. (5) is approximated as,

$$\int_0^{x_k} K(x_k, \tau) \varphi(\tau) d\tau = \sum_{j=0}^{k-1} \int_{x_j}^{x_{j+1}} K(x_k, \tau) \varphi(\tau) d\tau. \tag{8}$$

The function $\varphi(\tau)$ is approximated over the interval $[x_j, x_{j+1}]$ using the following formula $\varphi_{1j}(\tau) \approx \frac{(x_{j+1} - \tau)}{h} \varphi_j + \frac{(\tau - x_j)}{h} \varphi_{j+1}$, then

$$\int_0^{x_k} K(x_k, \tau) \varphi(\tau) d\tau = \sum_{j=0}^{k} U(k,j) \varphi_j, \tag{9}$$

where,

$$U(k,j) = \begin{cases} S(k,0) & , \quad j = 0 \\ S(k,j) + T(k, j-1) & , 1 \le j \le k \\ T(k, k-1) & , \quad j = k \end{cases} \tag{10}$$

$$S(k,j) = h \int_0^1 (1-p) K(kh, hp + jh) dp, \tag{11}$$

$$T(k,j) = h \int_0^1 p \, K(kh, hp + jh) dp, \tag{12}$$

Using Eq.(6) and Eq. (9), Eq. (5) is equivalent to the linear system of equations,

$$\sum_{j=0}^{k} \xi(k,j) \varphi_j = f_k \, , \quad k = 1, 2, \dots, n \, , \tag{13}$$

where,

$$\xi(k,j) = A(k,j) - U(k,j),$$

The approximated solution of the Eq. (4) can be obtained by solving the linear system given by Eq. (13). Now we discuss the Quadratic scheme (S2) as follows.

## 3.2 The Quadratic Scheme (S2)

In this scheme, the unknown function is approximated by linear interpolation function for $j = 0$ in the interval $[x_0, x_1]$ and for, $\ge 2$, the unknown function is approximated as quadratic interpolation function in the interval $[x_{j-1}, x_j]$. We follow the similar steps as described in the scheme S1 to derive the expression for the Quadratic scheme for Eq.(5). The derivations are obtained in parts and described below.

Part 1: Approximation of the left side of the Eq.(5) containing the Caputo derivative term:

We assume the initial condition, $D^\alpha \varphi(x_0) = 0$,

$$\frac{1}{\Gamma(1-\alpha)} \int_0^{x_k} (x_k - \tau)^{-\alpha} \varphi'(\tau) d\tau = \frac{1}{\Gamma(1-\alpha)} \sum_{j=1}^{k} \int_{x_{j-1}}^{x_j} (x_k - \tau)^{-\alpha} \varphi'(\tau) d\tau, \tag{14}$$

$$= \frac{1}{\Gamma(1-\alpha)} \int_{x_0}^{x_1} (x_k - \tau)^{-\alpha} \varphi'(\tau) d\tau + \frac{1}{\Gamma(1-\alpha)} \sum_{j=2}^{k} \int_{x_{j-1}}^{x_j} (x_k - \tau)^{-\alpha} \varphi'(\tau) d\tau. \tag{15}$$

On the first interval $[x_0, x_1]$, the linear interpolation is used for approximating the unknown function as like scheme S1 and for the other subintervals $(j \ge 2)$, the quadratic interpolation function $(\varphi_{2j-1}(\tau))'$ for the three points $(x_{j-2}, \varphi_{j-2}), (x_{j-1}, \varphi_{j-1}), (x_j, \varphi_j)$ is applied such that,

$$(\varphi_{2j-1}(\tau))' = \left( \frac{(\tau - x_{j-1})(\tau - x_j)}{2h^2} \varphi_{j-2} - \frac{(\tau - x_{j-2})(\tau - x_j)}{h^2} \varphi_{j-1} + \frac{(\tau - x_{j-2})(\tau - x_{j-1})}{2h^2} \varphi_j \right)'.$$

$$(\varphi_{2j-1}(\tau))' = \left( \frac{(2\tau - x_{j-1} - x_j)}{2h^2} \varphi_{j-2} - \frac{(2\tau - x_{j-2} - x_j)}{h^2} \varphi_{j-1} + \frac{(2\tau - x_{j-2} - x_{j-1})}{2h^2} \varphi_j \right). \tag{16}$$

Thus from Eq.(15) and Eq.(16), we have,

$$\approx \frac{1}{\Gamma(1-\alpha)} \int_{x_0}^{x_1} (x_k - \tau)^{-\alpha} \frac{\varphi_1 - \varphi_0}{h} d\tau + \frac{1}{\Gamma(1-\alpha)} \sum_{j=2}^{k} \int_{x_{j-1}}^{x_j} (x_k - \tau)^{-\alpha} \left( \varphi_{2j-1}(\tau) \right)' d\tau,$$

$$= \frac{1}{\Gamma(1-\alpha)} \int_{x_0}^{x_1} (x_k - \tau)^{-\alpha} \frac{\varphi_1 - \varphi_0}{h} d\tau + \frac{h^{-\alpha}}{\Gamma(2-\alpha)} \sum_{j=2}^{k} B(k-j)\, \varphi_{j-2} + C(k-j)\varphi_{j-1} + D(k-j)\varphi_j,$$

$$= \frac{h^{-\alpha} a_{k-1}}{\Gamma(2-\alpha)} (\varphi_1 - \varphi_0) + \frac{h^{-\alpha}}{\Gamma(2-\alpha)} \sum_{j=2}^{k} B(k-j)\, \varphi_{j-2} + C(k-j)\varphi_{j-1} + D(k-j)\varphi_j, \tag{17}$$

where,

$$a_{k-1} = k^{1-\alpha} - (k-1)^{1-\alpha},$$

$$B(k-j) = \frac{1}{(2-\alpha)} \left[ (k-j+1)^{1-\alpha} \left( k-j+\frac{\alpha}{2} \right) - (k-j)^{1-\alpha} \left( k-j-\frac{\alpha}{2}+1 \right) \right],$$

$$C(k-j) = \frac{2}{(2-\alpha)} \left[ (k-j)^{1-\alpha} (k-j-\alpha+2) - (k-j+1)^{2-\alpha} \right]$$

$$D(k-j) = \frac{1}{(2-\alpha)} \left[ (k-j+1)^{1-\alpha} \left( k-j-\frac{\alpha}{2}+2 \right) - (k-j)^{1-\alpha} \left( k-j-\frac{3\alpha}{2}+3 \right) \right].$$

Eq.(17) can be written as,

$$\frac{h^{-\alpha} a_{k-1}}{\Gamma(2-\alpha)} (\varphi_1 - \varphi_0) + \frac{h^{-\alpha}}{\Gamma(2-\alpha)} \sum_{j=2}^{k} B(k-j)\, \varphi_{j-2} + C(k-j)\varphi_{j-1} + D(k-j)\varphi_j = \frac{h^{-\alpha}}{\Gamma(2-\alpha)} \sum_{j=0}^{k} s_j\, \varphi_{k-j}.$$

Or, $\quad \frac{1}{\Gamma(1-\alpha)} \int_0^{x_k} (x_k - \tau)^{-\alpha} \varphi'(\tau) d\tau = \frac{h^{-\alpha}}{\Gamma(2-\alpha)} \sum_{j=0}^{k} s_j\, \varphi_{k-j}.$ \hfill (18)

Here the coefficients $s_j$ in Eq.(18) for different $k$ can be expressed as,

For $k=1$, $\quad \begin{cases} s_0 = a_0 \\ s_1 = -a_0 \end{cases}.$

For $k=2$, $\quad \begin{cases} s_0 = D(0) \\ s_1 = C(0) + a_1 \\ s_2 = B(0) - a_1 \end{cases}.$

For $k=3$, $\quad \begin{cases} s_0 = D(0) \\ s_1 = C(0) + D(1) \\ s_2 = B(0) + C(1) + a_2 \\ s_3 = B(1) - a_2 \end{cases}.$

And for $k \geq 4$ ,the coefficients have the following relations:

$$\begin{cases} s_0 = D(0) \\ s_1 = C(0) + D(1) \\ s_j = B(j-2) + C(j-1) + D(j)(2 \leq j \leq k-2). \\ s_{k-1} = B(k-3) + C(k-2) + a_{k-1} \\ s_k = B(k-2) - a_{k-1} \end{cases}$$

Part 2: Approximation of the integration term on the right side of Eq.(5):

The integration term on the right side of Eq.(5) is expressed as,

$$\int_0^{x_k} K(x_k, \tau) \varphi(\tau) d\tau = \sum_{j=1}^k \int_{x_{j-1}}^{x_j} K(x_k, \tau) \varphi(\tau) d\tau , \tag{19}$$

$$= \int_{x_0}^{x_1} K(x_k, \tau) \varphi(\tau) d\tau + \sum_{j=2}^k \int_{x_{j-1}}^{x_j} K(x_k, \tau) \varphi(\tau) d\tau. \tag{20}$$

In Eq.(20), the function $\varphi(\tau)$ is approximated over the interval $[x_0, x_1]$ using the following formula:

$$\varphi_{10}(\tau) \approx \frac{(x_1 - \tau)}{h} \varphi_0 + \frac{(\tau - x_0)}{h} \varphi_1 ,$$

$$\int_{x_0}^{x_1} K(x_k, \tau) \varphi(\tau) d\tau = a_k \varphi_0 + b_k \varphi_1 ,$$

where $a_k = h \int_0^1 K(kh, ph)(1-p) \, dp$ and $\quad b_k = h \int_0^1 K(kh, ph)p \, dp.$ \tag{21}

And for the $(j \geq 2)$ function $\varphi(\tau)$ is approximated over the interval $[x_{j-1}, x_j]$ for three points $(x_{j-2}, \varphi_{j-2}), (x_{j-1}, \varphi_{j-1}), (x_j, \varphi_j)$, we use the quadratic interpolation function $\varphi(\tau)$ as discussed in the discretization of the Part 1 such that,

$$\varphi_{2j-1}(\tau) = \left( \frac{(\tau - x_{j-1})(\tau - x_j)}{2h^2} \varphi_{j-2} - \frac{(\tau - x_{j-2})(\tau - x_j)}{h^2} \varphi_{j-1} + \frac{(\tau - x_{j-2})(\tau - x_{j-1})}{2h^2} \varphi_j \right).$$

Hence second integral of the Eq.(20) is expressed as,

$$\sum_{j=2}^k \int_{x_{j-1}}^{x_j} K(x_k, \tau) \varphi(\tau) d\tau \approx \sum_{j=2}^k M(k,j) \, \varphi_{j-2} + N(k,j) \varphi_{j-1} + O(k,j) \varphi_j. \tag{22}$$

Now using the transformation $\frac{\tau - x_{j-1}}{h} = p$, Eq.(20) takes form,

$$\sum_{j=1}^k \int_{x_{j-1}}^{x_j} K(x_k, \tau) \varphi(\tau) d\tau,$$

$$= a_k \varphi_0 + b_k \varphi_1 + \sum_{j=2}^k M(k,j) \, \varphi_{j-2} + N(k,j) \varphi_{j-1} + O(k,j) \varphi_j, \text{ where,}$$

$$M(k,j) = \frac{h}{2} \int_0^1 K(kh, hp + jh - h)\, p(p-1)\, dp,$$

$$N(k,j) = h \int_0^1 K(kh, hp + jh - h)\, (1 - p^2)\, dp,$$

$$O(k,j) = \frac{h}{2} \int_0^1 K(kh, hp + jh - h)\, p(p+1)\, dp,$$

Or,

$$\int_0^{x_k} K(x_k, \tau) \varphi(\tau)\, d\tau = \sum_{j=0}^k v_j\, \varphi_j. \tag{23}$$

The coefficients $v_j$ for different $k$ can be calculated as,

For $k = 1$, $\begin{cases} v_0 = a_1 \\ v_1 = b_1 \end{cases}$.

For $k = 2$, $\begin{cases} v_0 = M(2,2) + a_2 \\ v_1 = N(2,2) + b_2 \\ v_2 = O(2,2) \end{cases}$.

For $k = 3$, $\begin{cases} v_0 = M(3,2) + a_3 \\ v_1 = M(3,3) + N(3,2) + b_3 \\ v_2 = N(3,3) + O(3,2) \\ v_3 = O(3,3) \end{cases}$.

And for $k \geq 4$, the coefficients have the following relations:

$$\begin{cases} v_0 = M(k,2) + a_k \\ v_1 = M(k,3) + N(k,2) + b_k \\ v_j = M(k, j+2) + N(k, j+1) + O(k,j)\, (2 \leq j \leq k-2). \\ v_{k-1} = N(k,k) + O(k, k-1) \\ v_k = O(k,k) \end{cases}$$

Using Eq.(18) and Eq.(23), the Eq.(5) can be expressed as,

$$A_k(\varphi_0, \varphi_1, \varphi_2 \ldots \varphi_k) - C_k(\varphi_0, \varphi_1, \varphi_2 \ldots \varphi_k) = f_k, \quad k = 1, 2, \ldots, n, \tag{24}$$

where,

$$A_k(\varphi_0, \varphi_1, \varphi_2 \ldots \varphi_k) = \frac{h^{-\alpha}}{\Gamma(2-\alpha)} \sum_{j=0}^k s_j\, \varphi_{k-j},$$

$$C_k(\varphi_0, \varphi_1, \varphi_2 \ldots \varphi_k) = \sum_{j=0}^k v_j\, \varphi_j.$$

Now by solving the above system of linear equations given by Eq. (24), one obtains the desired approximate solution of the fractional integro-differential equation defined by Eq. (4).

### 3.3 The Quadratic-linear Scheme (S3)

Here, we discuss a hybrid scheme combining quadratic interpolation and linear interpolation approximation together to approximate the desired solution of Eq.(5). The left side of Eq. (5) (the term having Caputo derivative) is discretized using quadratic interpolation however the integration term on the right side of Eq.(5) is discretized using the linear interpolation approximation for the unknown function. We use the discretization of both the parts as discussed in schemes S1 and S2. Now from Eq. (18) and Eq. (9), Eq. (5) takes the form,

$$\frac{h^{-\alpha}}{\Gamma(2-\alpha)}\sum_{j=0}^{k} s_j \, \varphi_{k-j} - \sum_{j=0}^{k} U(k,j) \, \varphi_j = f_k \text{ , where, } k = 1,2,\dots,n. \tag{25}$$

Or,

$$A_k(\varphi_0 \, , \varphi_1, \varphi_2 \dots \varphi_k) - B_k(\varphi_0 \, , \varphi_1, \varphi_2 \dots \varphi_k) = f_k \text{ ,} \tag{26}$$

where,

$$A_k(\varphi_0 \, , \varphi_1, \varphi_2 \dots \varphi_k) = \frac{h^{-\alpha}}{\Gamma(2-\alpha)}\sum_{j=0}^{k} s_j \, \varphi_{k-j}, \tag{27}$$

$$B_k(\varphi_0 \, , \varphi_1, \varphi_2 \dots \varphi_k) = \sum_{j=0}^{k} U(k,j) \, \varphi_j \, . \tag{28}$$

By solving the above system given by Eq.(25), the approximate solution is calculated.

## 4. Error Estimate of the Approximations

In this section, the error estimates of the schemes S1, S2 and S3 are discussed as follows.

**Theorem 4.1 (S1)**: Suppose $\varphi(\tau) \in C^2[0 \, , x_k]$ and for any $\alpha( \, 0 < \alpha < 1)$. Suppose, $K(x,\tau), 0 \le x, \tau \le 1$ is continuous and $|K(x,\tau)| \le M$ where $M > 0$ . Let $E^k(\varphi, h, \alpha)$ be the error of the approximation to Eq.(4) using Linear scheme (S1) then,

$$\left|E^k(\varphi, h, \alpha)\right| \le \frac{1}{\Gamma(1-\alpha)}\left[\frac{1}{8} + \frac{\alpha}{(1-\alpha)(2-\alpha)}\right] \max_{x_0 \le \tau \le x_k} \left|\varphi''(\tau)\right| h^{2-\alpha} + \frac{M}{8} \max_{x_0 \le \tau \le x_k} \left|\varphi''(\tau)\right| x_k h^2.$$

**Proof**: From Eq.(4) and Eq.(13), we have,

$$\left|E^k(\varphi, \alpha)\right| \le \frac{1}{\Gamma(1-\alpha)}\left|\sum_{j=0}^{k-1} \int_{x_j}^{x_{j+1}}(x_k - \tau)^{-\alpha}\left(\varphi'(\tau) - \frac{\varphi_{j+1}-\varphi_j}{h}\right)d\tau\right| +$$

$$\left| \sum_{j=0}^{k-1} \int_{x_j}^{x_{j+1}} K(x_k,\tau) \left( \varphi(\tau) - p_{1j}(\tau) \right) d\tau \right|, \tag{29}$$

For simplicity, we use notation,

$$\left| E_1{}^k(\varphi,h,\alpha) \right| = \frac{1}{\Gamma(1-\alpha)} \left| \sum_{j=0}^{k-1} \int_{x_j}^{x_{j+1}} (x_k-\tau)^{-\alpha} \left( \varphi'(\tau) - \frac{\varphi_{j+1}-\varphi_j}{h} \right) d\tau \right|, \tag{30}$$

$$\left| E_2{}^k(\varphi,h,\alpha) \right| = \left| \sum_{j=0}^{k-1} \int_{x_j}^{x_{j+1}} K(x_k,\tau) \left( \varphi(\tau) - p_{1j}(\tau) \right) d\tau \right|. \tag{31}$$

Now we consider the first part given by Eq. (30),

$$\left| E_1{}^k(\varphi,h,\alpha) \right| = \frac{1}{\Gamma(1-\alpha)} \left| \sum_{j=0}^{k-1} \int_{x_j}^{x_{j+1}} (x_k-\tau)^{-\alpha} \left( \varphi'(\tau) - \frac{\varphi_{j+1}-\varphi_j}{h} \right) d\tau \right|,$$

$$= \frac{1}{\Gamma(1-\alpha)} \left| \sum_{j=0}^{k-1} \int_{x_j}^{x_{j+1}} (x_k-\tau)^{-\alpha} \left( \varphi(\tau) - \varphi_{1j}(\tau) \right)' d\tau \right|$$

$$= \frac{1}{\Gamma(1-\alpha)} \left| \sum_{j=0}^{k-1} \int_{x_j}^{x_{j+1}} (x_k-\tau)^{-\alpha} d\left( \varphi(\tau) - \varphi_{1j}(\tau) \right) \right|$$

$$= \frac{1}{\Gamma(1-\alpha)} \left| \sum_{j=0}^{k-1} \left( \varphi(\tau) - \varphi_{1j}(\tau) \right) (x_k-\tau)^{-\alpha} \big|_{\tau=x_j}^{x_{j+1}} - \frac{\alpha}{\Gamma(1-\alpha)} \sum_{j=0}^{k-1} \int_{x_j}^{x_{j+1}} (x_k-\tau)^{-\alpha-1} \left( \varphi(\tau) - \varphi_{1j}(\tau) \right) d\tau \right|$$

$$= \frac{\alpha}{\Gamma(1-\alpha)} \left| \sum_{j=0}^{k-1} \int_{x_j}^{x_{j+1}} (x_k-\tau)^{-\alpha-1} \left( \varphi(\tau) - \varphi_{1j}(\tau) \right) d\tau \right|$$

$$\leq \frac{\alpha}{\Gamma(1-\alpha)} \left| \sum_{j=0}^{k-2} \int_{x_j}^{x_{j+1}} (x_k-\tau)^{-\alpha-1} \left( \varphi(\tau) - \varphi_{1j}(\tau) \right) d\tau \right|$$

$$\qquad\qquad + \frac{\alpha}{\Gamma(1-\alpha)} \left| \int_{x_{k-1}}^{x_k} (x_k-\tau)^{-\alpha-1} \left( \varphi(\tau) - \varphi_{1k-1}(\tau) \right) d\tau \right|,$$

$$\leq \frac{1}{2} |\varphi''(\eta)| \frac{\alpha}{\Gamma(1-\alpha)} \left| \sum_{j=0}^{k-2} \int_{x_j}^{x_{j+1}} (x_k-\tau)^{-\alpha-1} (\tau-x_j)(\tau-x_{j+1}) d\tau \right|$$

$$\qquad\qquad + \frac{\alpha}{\Gamma(1-\alpha)} \frac{1}{2} |\varphi''(\xi)| \left| \int_{x_{k-1}}^{x_k} (x_k-\tau)^{-\alpha} (\tau-x_{k-1}) d\tau \right|$$

$$\leq \frac{1}{8} |\varphi''(\eta)| \frac{1}{\Gamma(1-\alpha)} h^{2-\alpha} + \frac{\alpha}{\Gamma(1-\alpha)} \frac{1}{2} |\varphi''(\xi)| \frac{h^{2-\alpha}}{(1-\alpha)(2-\alpha)}, \tag{32}$$

where, $\eta \in (x_0, x_{k-1})$, $\xi \in (x_{k-1}, x_k)$.

Now we take second part given by Eq. (31),

$$\left| E_2{}^k(\varphi,h,\alpha) \right| = \left| \sum_{j=0}^{k-1} \int_{x_j}^{x_{j+1}} K(x_k,\tau) \left( \varphi(\tau) - p_{1j}(\tau) \right) d\tau \right|,$$

$$\leq \frac{|\varphi''(\xi_1)|}{2} \left| \sum_{j=0}^{k-1} \int_{x_j}^{x_{j+1}} K(x_k,\tau)(\tau-x_j)(\tau-x_{j+1}) d\tau \right| \leq \frac{M}{8} \max_{x_0 \leq \tau \leq x_k} |\varphi''(\tau)| x_k h^2, \tag{33}$$

where, $\xi_1 \in (x_0, x_k)$.

From Eqs.(32-33) and Eq. (29),

$$|E^k(\varphi, h, \alpha)| \leq \frac{1}{\Gamma(1-\alpha)} \left[\frac{1}{8} + \frac{\alpha}{(1-\alpha)(2-\alpha)}\right] \max_{x_0 \leq \tau \leq x_k} |\varphi''(\tau)| h^{2-\alpha} + \frac{M}{8} \max_{x_0 \leq x \leq x_k} |\varphi''(\tau)| x_k h^2.$$

This completes the proof.

**Theorem 4.2 (S2):** Suppose $\varphi(\tau) \in C^3[0, x_k]$ and for any $\alpha$ ( $0 < \alpha < 1$). Suppose, $K(x, \tau), 0 \leq x, \tau \leq 1$ is continuous and $|K(x, \tau)| \leq M$ where $M > 0$. Let $G^k(\varphi, h, \alpha)$ be the error of the approximation to Eq. (4) using Quadratic scheme (S2) then,

(i) For, $= 1$ , $|G^1(\varphi, h, \alpha)| \leq \frac{\alpha}{2\Gamma(3-\alpha)} \max_{x_0 \leq \tau \leq x_1} |\varphi''(\tau)| h^{2-\alpha} + \frac{M}{2} \max_{x_0 \leq \tau \leq x_1} |\varphi''(\tau)| x_1 h^2$, and

(ii) For, $k \geq 2$ $\quad |G^k(\varphi, h, \alpha)| \leq$

$$\frac{1}{\Gamma(1-\alpha)} \left(\frac{\alpha}{12} \max_{x_0 \leq \tau \leq x_k} |\varphi''(\tau)| (x_k - x_1)^{-\alpha-1} h^3 + \left[\frac{1}{12} + \frac{\alpha}{3(1-\alpha)(2-\alpha)} \left(\frac{1}{2} + \frac{1}{3-\alpha}\right)\right] \max_{x_0 \leq \tau \leq x_k} |\varphi'''(\tau)| h^{3-\alpha}\right)$$
$$+ \frac{M}{2} \max_{x_0 \leq \tau \leq x_1} |\varphi''(\tau)| x_1 h^2 + \frac{M}{12} \max_{x_0 \leq \tau \leq x_k} |\varphi'''(\tau)| (x_k - x_1) h^3.$$

**Proof:** Part (i) From Eq.(4) and Eq.(24), we have the error of the approximation for $k = 1$,

$$|G^1(\varphi, h, \alpha)| \leq \frac{1}{\Gamma(1-\alpha)} \left|\int_{x_0}^{x_1} (x_1 - \tau)^{-\alpha} \left(\varphi'(\tau) - \frac{\varphi_1 - \varphi_0}{h}\right) d\tau\right| + \left|\int_{x_0}^{x_1} K(x_1, \tau) \left(\varphi(\tau) - p_{1j}(\tau)\right) d\tau\right|. \quad (34)$$

For simplicity, we introduce the notations for each right side part of Eq. (34) as,

$$|G_1^1(\varphi, h, \alpha)| = \frac{1}{\Gamma(1-\alpha)} \left|\int_{x_0}^{x_1} (x_1 - \tau)^{-\alpha} \left(\varphi'(\tau) - \frac{\varphi_1 - \varphi_0}{h}\right) d\tau\right|, \quad (35)$$

$$|G_2^1(\varphi, h, \alpha)| = \left|\int_{x_0}^{x_1} K(x_1, \tau) \left(\varphi(\tau) - p_{1j}(\tau)\right) d\tau\right|. \quad (36)$$

Now we consider the first part given by Eq. (35),

$$|G_1^1(\varphi, h, \alpha)| = \frac{1}{\Gamma(1-\alpha)} \left|\int_{x_0}^{x_1} (x_1 - \tau)^{-\alpha} \left(\varphi'(\tau) - \frac{\varphi_1 - \varphi_0}{h}\right) d\tau\right| \leq \frac{\alpha}{2\Gamma(3-\alpha)} \max_{x_0 \leq \tau \leq x_1} |\varphi''(\tau)| h^{2-\alpha}, \quad (37)$$

$$|G_2^1(\varphi, h, \alpha)| = \left|\int_{x_0}^{x_1} K(x_1, \tau) \left(\varphi(\tau) - p_{1j}(\tau)\right) d\tau\right| \leq \frac{M}{2} \max_{x_0 \leq \tau \leq x_1} |\varphi''(\tau)| x_1 h^2. \quad (38)$$

From Eqs. (37)-(38) and Eq. (34),

$$|G^1(\varphi, h, \alpha)| \leq \frac{\alpha}{2\Gamma(3-\alpha)} \max_{x_0 \leq \tau \leq x_1} |\varphi''(\tau)| h^{2-\alpha} + \frac{M}{2} \max_{x_0 \leq \tau \leq x_1} |\varphi''(\tau)| x_1 h^2.$$

This completes the proof for Part (i). The proof of Part (ii) is given below.

Part (ii) From Eq.(4) and Eq.(24), we have, for $k \geq 2$

$$\left| G^k(\varphi, h, \alpha) \right| \leq$$

$$\frac{1}{\Gamma(1-\alpha)} \left( \left| \int_{x_0}^{x_1} (x_k - \tau)^{-\alpha} \left( \varphi(\tau) - \varphi_{10}(\tau) \right)' d\tau \right| + \left| \sum_{j=2}^{k} \int_{x_{j-1}}^{x_j} (x_k - \tau)^{-\alpha} \left( \varphi(\tau) - \varphi_{2j-1}(\tau) \right)' d\tau \right| \right)$$

$$+ \left| \int_{x_0}^{x_1} K(x_k, \tau) \left( \varphi(\tau) - p_{1j}(\tau) \right) d\tau \right| + \left| \sum_{j=2}^{k} \int_{x_{j-1}}^{x_j} K(x_k, \tau) \left( \varphi(\tau) - p_{1j}(\tau) \right) d\tau \right|. \tag{39}$$

We denote Eq. (39) as,

$$\left| G^k(\varphi, h, \alpha) \right| \leq \left| G_1{}^k(\varphi, h, \alpha) \right| + \left| G_2{}^k(\varphi, h, \alpha) \right|, \tag{40}$$

where,

$$\left| G_1{}^k(\varphi, h, \alpha) \right| = \frac{1}{\Gamma(1-\alpha)} \left( \left| \int_{x_0}^{x_1} (x_k - \tau)^{-\alpha} \left( \varphi(\tau) - \varphi_{10}(\tau) \right)' d\tau \right| + \left| \sum_{j=2}^{k} \int_{x_{j-1}}^{x_j} \frac{\left( \varphi(\tau) - \varphi_{2j-1}(\tau) \right)'}{(x_k - \tau)^{\alpha}} d\tau \right| \right) \tag{41}$$

$$\left| G_2{}^k(\varphi, h, \alpha) \right| = \left| \int_{x_0}^{x_1} K(x_k, \tau) \left( \varphi(\tau) - \varphi_{10}(\tau) \right) d\tau \right| + \left| \sum_{j=2}^{k} \int_{x_{j-1}}^{x_j} K(x_k, \tau) \left( \varphi(\tau) - \varphi_{2j-1}(\tau) \right) d\tau \right| \tag{42}$$

Now we consider Eq. (41),

$$\left| G_1{}^k(\varphi, h, \alpha) \right| = \frac{1}{\Gamma(1-\alpha)} \left( \left| \int_{x_0}^{x_1} (x_k - \tau)^{-\alpha} \left( \varphi(\tau) - \varphi_{10}(\tau) \right)' d\tau \right| + \left| \sum_{j=2}^{k} \int_{x_{j-1}}^{x_j} \frac{\left( \varphi(\tau) - \varphi_{2j-1}(\tau) \right)'}{(x_k - \tau)^{\alpha}} d\tau \right| \right),$$

$$\leq \frac{\alpha}{\Gamma(1-\alpha)} \left[ \left| \int_{x_0}^{x_1} \left( \varphi(\tau) - \varphi_{10}(\tau) \right) (x_k - \tau)^{-\alpha-1} d\tau \right| + \left| \sum_{j=2}^{k} \int_{x_{j-1}}^{x_j} \left( \varphi(\tau) - \varphi_{2j-1}(\tau) \right) (x_k - \tau)^{-\alpha-1} d\tau \right| \right]. \tag{43}$$

Since,

$$\left| \int_{x_0}^{x_1} \left( \varphi(\tau) - \varphi_{10}(\tau) \right) (x_k - \tau)^{-\alpha-1} d\tau \right|$$

$$= \left| \int_{x_0}^{x_1} \frac{\varphi''(\xi_2)}{2} (\tau - x_0)(\tau - x_1)(x_k - \tau)^{-\alpha-1} d\tau \right|,$$

$$= \left| \frac{\varphi''(\eta_1)}{2} \int_{x_0}^{x_1} (\tau - x_0)(x_1 - \tau)(x_k - \tau)^{-\alpha-1} d\tau \right|,$$

$$\leq \frac{1}{12}|\varphi''(\eta_1)|(x_k - x_1)^{-\alpha-1}h^3, \tag{44}$$

where $\eta_1 \in (x_0, x_1)$. We also know that

$$\left|\sum_{j=2}^{k-1} \int_{x_{j-1}}^{x_j} \left(\varphi(\tau) - \varphi_{2j-1}(\tau)\right)(x_k - \tau)^{-\alpha-1}d\tau\right|,$$

$$= \left|\sum_{j=2}^{k-1} \int_{x_{j-1}}^{x_j} \frac{\varphi'''(\vartheta_j)}{6}(\tau - x_{j-2})(\tau - x_{j-1})(\tau - x_j)(x_k - \tau)^{-\alpha-1}d\tau\right|,$$

$$= \frac{1}{6}\left|\sum_{j=2}^{k-1} \varphi'''(\eta_j) \int_{x_{j-1}}^{x_j} (\tau - x_{j-2})(\tau - x_{j-1})(x_j - \tau)(x_k - \tau)^{-\alpha-1}d\tau\right|,$$

$$= \frac{1}{6}|\varphi'''(\eta)| \sum_{j=2}^{k-1} \int_{x_{j-1}}^{x_j} (\tau - x_{j-2})(\tau - x_{j-1})(x_j - \tau)(x_k - \tau)^{-\alpha-1}d\tau,$$

$$\leq \frac{1}{12}|\varphi'''(\eta)|h^3 \int_{x_1}^{x_{k-1}} (x_k - \tau)^{-\alpha-1}d\tau,$$

$$\leq \frac{1}{12}|\varphi'''(\eta)|h^{3-\alpha}, \tag{45}$$

where $\eta_j \in (x_{j-2}, x_j), 2 \leq j \leq k-1, \eta \in (x_0, x_{k-1})$, and

$$\int_{x_{k-1}}^{x_k} \left(\varphi(\tau) - \varphi_{2k-1}(\tau)\right)(x_k - \tau)^{-\alpha-1}d\tau,$$

$$= \int_{x_{k-1}}^{x_k} \frac{\varphi'''(\vartheta_k)}{6}(\tau - x_{k-2})(\tau - x_{k-1})(\tau - x_k)(x_k - \tau)^{-\alpha-1}d\tau,$$

$$= -\frac{\varphi'''(\eta_k)}{6} \int_{x_{k-1}}^{x_k} (\tau - x_{k-2})(\tau - x_{k-1})(x_k - \tau)^{-\alpha}d\tau,$$

$$= -\frac{1}{3}\frac{1}{(2-\alpha)(1-\alpha)}\left(\frac{1}{2} + \frac{1}{3-\alpha}\right)\varphi'''(\eta_k)h^{3-\alpha}. \tag{46}$$

From Eqs. (42)-(46) and (40) we get,

$$\left|G_1{}^k(\varphi, h, \alpha)\right| \leq$$

$$\left(\frac{\alpha}{12} \max_{x_0 \leq x \leq x_1}|\varphi''(\tau)|(x_k - x_1)^{-\alpha-1}h^3 + \left[\frac{1}{12} + \frac{\alpha}{3(1-\alpha)(2-\alpha)}\left(\frac{1}{2} + \frac{1}{3-\alpha}\right)\right] \max_{x_0 \leq \tau \leq x_k} \frac{1}{\Gamma(1-\alpha)}|\varphi'''(\tau)|h^{3-\alpha}\right). \tag{47}$$

Now we consider Eq. (42),

$$\left|G_2{}^k(\varphi, h, \alpha)\right| = \left|\int_{x_0}^{x_1} K(x_k, \tau)\left(\varphi - \varphi_{10}(\tau)\right)d\tau\right| + \left|\sum_{j=2}^k \int_{x_{j-1}}^{x_j} K(x_k, \tau)\left(\varphi - \varphi_{2j-1}(\tau)\right)d\tau\right|$$

$$= \frac{M}{2} \max_{x_0 \leq \tau \leq x_1}|\varphi''(\tau)|x_1 h^2 + \frac{1}{6} \max_{x_0 \leq \tau \leq x_k}|\varphi'''(\tau)|\left|\sum_{j=2}^k \int_{x_{j-1}}^{x_j} K(x_k, \tau)(\tau - x_{j-2})(\tau - x_{j-1})(\tau - x_j)d\tau\right|,$$

$$\leq \frac{M}{2} \max_{x_0 \leq \tau \leq x_1} |\varphi''(\tau)| x_1 h^2 + \frac{M}{12} \max_{x_0 \leq \tau \leq x_k} |\varphi'''(\tau)| (x_k - x_1) h^3 \ . \tag{48}$$

From Eqs.(47) , Eq. (48) and Eq.(40),

$$|G^k(\varphi, h, \alpha)| \leq$$

$$\frac{1}{\Gamma(1-\alpha)} \left( \frac{\alpha}{12} \max_{x_0 \leq \tau \leq x_1} |\varphi''(\tau)| (x_k - x_1)^{-\alpha-1} h^3 + \left[ \frac{1}{12} + \frac{\alpha}{3(1-\alpha)(2-\alpha)} \left( \frac{1}{2} + \frac{1}{3-\alpha} \right) \right] \max_{x_0 \leq \tau \leq x_k} |\varphi'''(\tau)| h^{3-\alpha} \right)$$

$$+ \frac{M}{2} \max_{x_0 \leq \tau \leq x_1} |\varphi''(\tau)| x_1 h^2 + \frac{M}{12} \max_{x_0 \leq \tau \leq x_k} |\varphi'''(\tau)| (x_k - x_1) h^3.$$

This completes the proof.

**Theorem 4.3 (S3)**: Suppose $\varphi(\tau) \in C^3[0, x_k]$ and for any $\alpha( 0 < \alpha < 1)$. Suppose, $K(x, \tau), 0 \leq x, \tau \leq 1$ is continuous and $|K(x, \tau)| \leq M$ where $M > 0$ . Let $F^k(\varphi, h, \alpha)$ be the error of the approximation to Eq.(4) using Quadratic-linear scheme (S3) then,

(i) For, $= 1$ , $|F^1(\varphi, h, \alpha)| \leq \frac{\alpha}{2\Gamma(3-\alpha)} \max_{x_0 \leq \tau \leq x_1} |\varphi''(\tau)| h^{2-\alpha} + \frac{M}{2} \max_{x_0 \leq \tau \leq x_1} |\varphi''(\tau)| x_1 h^2.$

(ii) For, $k \geq 2$, $\quad |F^k(\varphi, h, \alpha)| \leq$

$$\left( \frac{\alpha}{12} \max_{x_0 \leq \tau \leq x_1} |\varphi''(\tau)| (x_k - x_1)^{-\alpha-1} h^3 + \left[ \frac{1}{12} + \frac{\alpha}{3(1-\alpha)(2-\alpha)} \left( \frac{1}{2} + \frac{1}{3-\alpha} \right) \right] \max_{x_0 \leq \tau \leq x_k} \frac{1}{\Gamma(1-\alpha)} |\varphi'''(\tau)| h^{3-\alpha} \right)$$

$$+ \frac{M}{2} \max_{x_0 \leq \tau \leq x_k} |\varphi''(\tau)| x_k h^2.$$

**Proof:**

(i) From Eq.(4) and Eq.(25), we have,

$$|F^1(\varphi, h, \alpha)| \leq \frac{1}{\Gamma(1-\alpha)} \left| \int_{x_0}^{x_1} (x_1 - \tau)^{-\alpha} \left( \varphi'(\tau) - \frac{\varphi_1 - \varphi_0}{h} \right) d\tau \right| + \left| \int_{x_0}^{x_1} K(x_1, \tau) \left( \varphi(\tau) - \varphi_{10}(\tau) \right) d\tau \right| \tag{49}$$

For simplicity, we use notation,

$$|F_1{}^1(\varphi, h, \alpha)| = \frac{1}{\Gamma(1-\alpha)} \left| \int_{x_0}^{x_1} (x_1 - \tau)^{-\alpha} \left( \varphi'(\tau) - \frac{\varphi_1 - \varphi_0}{h} \right) d\tau \right|, \tag{50}$$

$$|F_2{}^1(\varphi, h, \alpha)| = \left| \int_{x_0}^{x_1} K(x_1, \tau) \left( \varphi(\tau) - \varphi_{10}(\tau) \right) d\tau \right|, \tag{51}$$

such that,

$$|F^1(\varphi, h, \alpha)| \leq |F_1{}^1(\varphi, h, \alpha)| + |F_2{}^1(\varphi, h, \alpha)|. \tag{52}$$

Proof of this theorem could be obtained using Theorem 4.2 and Theorem 4.1.

Now we consider Eq. (50),

$$\left|F_1{}^1(\varphi, h, \alpha)\right| = \frac{1}{\Gamma(1-\alpha)}\left|\int_{x_0}^{x_1}(x_1-\tau)^{-\alpha}\left(\varphi'(\tau) - \frac{\varphi_1-\varphi_0}{h}\right)d\tau\right| \leq \frac{\alpha}{2\Gamma(3-\alpha)}\max_{x_0\leq\tau\leq x_1}|\varphi''(\tau)|h^{2-\alpha}. \quad (53)$$

$$\left|F_2{}^1(\varphi, h, \alpha)\right| = \left|\int_{x_0}^{x_1}K(x_1,\tau)\big(\varphi(\tau)-\varphi_{10}(\tau)\big)d\tau\right| \leq \frac{M}{2}\max_{x_0\leq\tau\leq x_1}|\varphi''(\tau)|x_1 h^2. \quad (54)$$

From Eqs.(52)-(54),

$$|F^1(\varphi, h, \alpha)| \leq \frac{\alpha}{2\Gamma(3-\alpha)}\max_{x_0\leq\tau\leq x_1}|\varphi''(\tau)|h^{2-\alpha} + \frac{M}{2}\max_{x_0\leq\tau\leq x_1}|\varphi''(\tau)|x_1 h^2.$$

(ii) From Eq(4) and Eq.(25), we have, for $k \geq 2$

$$\left|F^k(\varphi, h, \alpha)\right|$$

$$\leq \frac{1}{\Gamma(1-\alpha)}\left(\left|\int_{x_0}^{x_1}(x_k-\tau)^{-\alpha}\big(\varphi(\tau)-\varphi_{10}(\tau)\big)'d\tau\right| + \left|\sum_{j=2}^{k}\int_{x_{j-1}}^{x_j}(x_k-\tau)^{-\alpha}\left(\varphi(\tau)-\varphi_{2j-1}(\tau)\right)'d\tau\right|\right)$$

$$\quad + \left|\sum_{j=0}^{k-1}\int_{x_j}^{x_{j+1}}K(x_k,\tau)\left(\varphi(\tau)-p_{1j}(\tau)\right)d\tau\right|. \quad (55)$$

Eq. (55) can be written directly using Eq.(33) and Eq.(41),

$$\left|F^k(\varphi, h, \alpha)\right| \leq$$

$$\left(\frac{\alpha}{12}\max_{x_0\leq\tau\leq x_1}|\varphi''(\tau)|(x_k-x_1)^{-\alpha-1}h^3 + \left[\frac{1}{12}+\frac{\alpha}{3(1-\alpha)(2-\alpha)}\left(\frac{1}{2}+\frac{1}{3-\alpha}\right)\right]\max_{x_0\leq\tau\leq x_k}\frac{1}{\Gamma(1-\alpha)}|\varphi'''(\tau)|h^{3-\alpha}\right)$$

$$\quad + \frac{M}{2}\max_{x_0\leq\tau\leq x_k}|\varphi''(\tau)|x_k h^2.$$

This completes the proof.

# 5 Numerical Results

Here, we take the test examples from the literature and investigate the performance of the presented schemes from section 4. The test examples are solved in [21]. The numerical solution obtained using the discussed schemes S1, S2 and S3 are presented through the tables. Further, the convergence order of the presented schemes is calculated varying the step size $h$. For calculating the convergence order (CO) in each case, the maximum absolute error (MAE) is calculated using the formula,

$$\text{E(x)} = \varphi_{exact}(x) - \varphi_{Numerical}(x),$$

$$MAE(h) := max\{|E(x)|, x \in \{x_0, x_1, \dots, x_k\}\},$$

and then the convergence order(CO) is obtained using the formula CO =lg[MAE(h)/MAE(h/2)]/lg(2).

**Test example 5.1:** Consider the following linear fractional integro-differential equation [21],

$$D^{1/2}\varphi(x) = \frac{(8/3)x^{3/2} - 2x^{1/2}}{\sqrt{\pi}} - \frac{3x^5 - 4x^4}{12} + \int_0^x x\tau\varphi(\tau)d\tau \, , \, , \quad 0 \le x, \tau \le 1,$$

subject to $\varphi(0) = 0$, having exact solution $\varphi(x) = x^2 - x$.

The numerical solutions of test example 5.1 using schemes S1, S2 and S3 are presented through Tables 1-2 for varying the step size $h = 1/5$ and $h = 1/10$ respectively. Further, maximum absolute errors and the convergence orders for test example 5.1 are provided in Tables 3-5 for the schemes S1, S2 and S3. From the Tables 3-5, it is observed that the scheme S3 which is the combination of the Quadratic and Linear schemes provide comparatively better convergence order than the schemes S1 and S2.

Table 1 Numerical solutions obtained using schemes S1, S2 and S3 for test example 5.1 for $n = 5$.

| $x_j$ | Exact solution | S1 | S2 | S3 |
|-------|----------------|-----------|-----------|-----------|
| 0.0 | 0.00 | 0.00 | 0.00 | 0.00 |
| 0.2 | -0.16 | -0.146642 | -0.146720 | -0.146642 |
| 0.4 | -0.24 | -0.217607 | -0.228706 | -0.228487 |
| 0.6 | -0.24 | -0.209663 | -0.231148 | -0.230553 |
| 0.8 | -0.16 | -0.120963 | -0.152223 | -0.150925 |
| 1.0 | 0.00 | 0.0505046 | 0.0078693 | 0.0104518 |

Table 2 Numerical solutions obtained using schemes S1, S2 and S3 for test example 1 for $n = 10$.

| $x_j$ | Exact solution | S1 | S2 | S3 |
|---|---|---|---|---|
| 0.0 | 0.00 | 0.00 | 0.00 | 0.00 |
| 0.2 | -0.16 | -0.154475 | -0.157177 | -0.157166 |
| 0.4 | -0.24 | -0.231231 | -0.238252 | -0.238209 |
| 0.6 | -0.24 | -0.228457 | -0.238611 | -0.23849 |
| 0.8 | -0.16 | -0.145454 | -0.158719 | -0.154838 |
| 1.0 | 0.00 | 0.0184683 | 0.00132716 | 0.00190422 |

Table 3 Maximum absolute error(MAE) and convergence order(CO) for test example 5.1 using scheme S1.

| $h$ | MAE | CO |
|---|---|---|
| $\frac{1}{5}$ | $5.05046 \times 10^{-2}$ | |
| $\frac{1}{10}$ | $1.84683 \times 10^{-2}$ | 1.45136 |
| $\frac{1}{20}$ | $6.70253 \times 10^{-3}$ | 1.46227 |
| $\frac{1}{40}$ | $2.41481 \times 10^{-3}$ | 1.4728 |
| $\frac{1}{80}$ | $8.65157 \times 10^{-4}$ | 1.48087 |

Table 4 Maximum absolute error(MAE) and convergence order(CO) for test example 5.1 using scheme S2.

| $h$ | MAE | CO |
|:---:|:---:|:---:|
| $\dfrac{1}{5}$ | $1.32804 \times 10^{-2}$ | |
| $\dfrac{1}{10}$ | $3.32983 \times 10^{-3}$ | 1.9958 |
| $\dfrac{1}{20}$ | $8.33153 \times 10^{-4}$ | 1.9988 |
| $\dfrac{1}{40}$ | $2.08325 \times 10^{-4}$ | 1.9998 |
| $\dfrac{1}{80}$ | $5.20829 \times 10^{-5}$ | 1.9995 |

Table 5 Maximum absolute error(MAE) and convergence order(CO) for test example 5.1 using scheme S3.

| $h$ | MAE | CO |
|:---:|:---:|:---:|
| $\dfrac{1}{5}$ | $1.3358 \times 10^{-2}$ | |
| $\dfrac{1}{10}$ | $3.33388 \times 10^{-3}$ | 2.00243 |
| $\dfrac{1}{20}$ | $8.33345 \times 10^{-4}$ | 2.00022 |
| $\dfrac{1}{40}$ | $2.08334 \times 10^{-4}$ | 2.00002 |
| $\dfrac{1}{80}$ | $5.20833 \times 10^{-5}$ | 2.00001 |

**Test example 5.2**: Here, we consider the fractional integro-differential equation [21],

$D^{5/6}\varphi(x) = f(x) + \int_0^x xe^\tau \varphi(\tau)d\tau,$, $\quad 0 \le x, \tau \le 1,$

subject to $\varphi(0) = 0$, where, $f(x) = -\frac{3}{91}\frac{\Gamma(5/6)x^{(1/6)}(-91+216x^2)}{\pi} + 5x - xe^x(5 - 5x + 3x^2 - x^3)$.

This has exact solution $\varphi(x) = x - x^3$.

The test example 5.2 is solved using the numerical schemes S1, S2 and S3 and the obtained approximate solutions are presented through Tables 6-7 varying the step size $h = 1/5$ and $h = 1/10$ respectively. However, Tables 8-10 present the maximum absolute errors and convergence orders for each schemes. One can observe that the scheme S1 achieves convergence order more than 1 and schemes S2 and S3 obtain the convergence order greater than 2. It is to mention here that the quadratic scheme S2 in this case achieves better convergence order (Table 9) than the schemes S3 (Table 10). This phenomena may be due to the behaviour of the kernel, $k(x,\tau) = xe^\tau$ in the right side part of the fractional integro-differential equation.

Table 6 Numerical solutions obtained using schemes S1, S2 and S3 for test example 5.2 for $n = 5$.

| $x_j$ | Exact solution | S1 | S2 | S3 |
|-------|---------------|-----|-----|-----|
| 0.0 | 0.00 | 0.00 | 0.00 | 0.00 |
| 0.2 | 0.192 | 0.180928 | 0.180928 | 0.180928 |
| 0.4 | 0.336 | 0.297013 | 0.314326 | 0.314249 |
| 0.6 | 0.384 | 0.298545 | 0.351121 | 0.350619 |
| 0.8 | 0.288 | 0.130392 | 0.241245 | 0.239352 |
| 1.0 | 0.000 | -0.273753 | -0.0670569 | -0.0726013 |

Table 7 Numerical solutions obtained using schemes S1, S2 and S3 for test example 5.2 for $n = 10$.

| $x_j$ | Exact solution | S1 | S2 | S3 |
|------|---------------|-----|-----|-----|
| 0.0 | 0.00 | 0.00 | 0.00 | 0.00 |
| 0.2 | 0.192 | 0.187244 | 0.189354 | 0.189353 |
| 0.4 | 0.336 | 0.318846 | 0.330992 | 0.330968 |
| 0.6 | 0.384 | 0.346531 | 0.376553 | 0.376414 |
| 0.8 | 0.288 | 0.219913 | 0.277566 | 0.277066 |
| 1.0 | 0.000 | -0.115363 | -0.0147491 | -0.0161695 |

Table 8 Maximum absolute error(MAE) and convergence order(CO) for test example 5.2 using scheme S1.

| $h$ | MAE | CO |
|-----|-----|-----|
| $\frac{1}{5}$ | $2.73753 \times 10^{-1}$ | |
| $\frac{1}{10}$ | $1.15363 \times 10^{-1}$ | 1.24669 |
| $\frac{1}{20}$ | $5.05187 \times 10^{-2}$ | 1.19129 |
| $\frac{1}{40}$ | $2.24081 \times 10^{-2}$ | 1.1728 |
| $\frac{1}{80}$ | $9.97983 \times 10^{-3}$ | 1.16693 |

Table 9 Maximum absolute error(MAE) and convergence order(CO) for test example 5.2 using scheme S2.

| $h$ | MAE | CO |
|---|---|---|
| $\frac{1}{5}$ | $6.70569 \times 10^{-2}$ | |
| $\frac{1}{10}$ | $1.47491 \times 10^{-2}$ | 2.18476 |
| $\frac{1}{20}$ | $3.28518 \times 10^{-3}$ | 2.16658 |
| $\frac{1}{40}$ | $7.33189 \times 10^{-4}$ | 2.16373 |
| $\frac{1}{80}$ | $1.63575 \times 10^{-4}$ | 2.16423 |

Table 10 Maximum absolute error(MAE) and convergence order(CO) for test example 5.2 using scheme S3.

| $h$ | MAE | CO |
|---|---|---|
| $\frac{1}{5}$ | $7.26013 \times 10^{-2}$ | |
| $\frac{1}{10}$ | $1.61695 \times 10^{-2}$ | 2.16672 |
| $\frac{1}{20}$ | $3.65192 \times 10^{-3}$ | 2.14655 |
| $\frac{1}{40}$ | $8.26938 \times 10^{-4}$ | 2.1428 |
| $\frac{1}{80}$ | $1.87313 \times 10^{-4}$ | 2.14233 |

**Test example 5.3:** In this example, we consider the fractional integro-differential equation similar to the one considered in [15] such that,

$$D^{1/3}\varphi(x) = \frac{3\sqrt{\pi}x^{7/6}}{4\Gamma(13/6)} - \frac{2}{63}x^{9/2}(9 + 7x^2) + \int_0^x (x\tau + x^2\tau^2)\varphi(\tau)d\tau, \quad 0 \le x, \tau \le 1,$$

subject to $\varphi(0) = 0$ with the exact solution $\varphi(x) = x^{3/2}$.

The numerical schemes S1, S2 and S3 are performed on this example and the obtained numerical results varying the step size $h$ are shown through Tables. Tables 11-12 present the obtained numerical solution using schemes S1, S2 and S3 for step sizes $h = 1/5$ and $h = 1/10$ respectively. And Tables 13-15 represent the maximum absolute errors and convergence orders for each schemes. From Tables 13-15, we observe that the scheme S1 performs comparatively better than the schemes S2 and S3. The occurrence of such behaviour may be due to the appearance of the fractional powers terms in the right side of the integro-differential equations and approximation of such fractional terms becomes difficult.

Table 11 Numerical solutions obtained using schemes S1, S2 and S3 for test example 5.3 for $n = 5$.

| $x_j$ | Exact solution | S1 | S2 | S3 |
|-------|---------------|------|------|------|
| 0.0 | 0.00 | 0.00 | 0.00 | 0.00 |
| 0.2 | 0.089443 | 0.099203 | 0.099203 | 0.099203 |
| 0.4 | 0.252982 | 0.263841 | 0.257785 | 0.257849 |
| 0.6 | 0.464758 | 0.476607 | 0.467637 | 0.467942 |
| 0.8 | 0.715542 | 0.72966 | 0.718179 | 0.719065 |
| 1.0 | 1.00 | 1.01967 | 1.00327 | 1.00547 |

Table 12 Numerical solutions obtained using schemes S1, S2 and S3 for test example 5.3 for $n = 10$.

| $x_j$ | Exact solution | S1 | S2 | S3 |
|---|---|---|---|---|
| 0.0 | 0.00 | 0.00 | 0.00 | 0.00 |
| 0.2 | 0.089443 | 0.093188 | 0.091088 | 0.09109 |
| 0.4 | 0.252982 | 0.256826 | 0.253615 | 0.253637 |
| 0.6 | 0.464758 | 0.46881 | 0.465257 | 0.465341 |
| 0.8 | 0.715542 | 0.72024 | 0.716039 | 0.716264 |
| 1.0 | 1.00 | 1.00634 | 1.00061 | 1.00115 |

Table 13 Maximum absolute error (MAE) and convergence order (CO) for test example 5.3 using scheme S1.

| $h$ | MAE | CO |
|---|---|---|
| $\dfrac{1}{5}$ | $1.96715 \times 10^{-2}$ | |
| $\dfrac{1}{10}$ | $6.34357 \times 10^{-3}$ | 1.63274 |
| $\dfrac{1}{20}$ | $2.05203 \times 10^{-3}$ | 1.62824 |
| $\dfrac{1}{40}$ | $6.61465 \times 10^{-4}$ | 1.63332 |
| $\dfrac{1}{80}$ | $2.12232 \times 10^{-4}$ | 1.64002 |

Table 14 Maximum absolute error (MAE) and convergence order (CO) for test example 5.3 using scheme S2.

| $h$ | MAE | CO |
|---|---|---|
| $\dfrac{1}{5}$ | $9.76046 \times 10^{-3}$ | |
| $\dfrac{1}{10}$ | $3.44045 \times 10^{-3}$ | 1.50435 |
| $\dfrac{1}{20}$ | $1.21603 \times 10^{-3}$ | 1.50042 |
| $\dfrac{1}{40}$ | $4.29918 \times 10^{-4}$ | 1.50005 |
| $\dfrac{1}{80}$ | $1.51999 \times 10^{-4}$ | 1.5 |

Table 15 Maximum absolute error (MAE) and convergence order (CO) for test example 5.3 using scheme S3.

| $h$ | MAE | CO |
|---|---|---|
| $\dfrac{1}{5}$ | $9.76046 \times 10^{-3}$ | |
| $\dfrac{1}{10}$ | $3.44045 \times 10^{-3}$ | 1.50435 |
| $\dfrac{1}{20}$ | $1.21603 \times 10^{-3}$ | 1.50042 |
| $\dfrac{1}{40}$ | $4.29918 \times 10^{-4}$ | 1.50005 |
| $\dfrac{1}{80}$ | $1.51999 \times 10^{-4}$ | 1.5 |

## Conclusions

We studied comparative study of the different approximations schemes such as Linear, Quadratic and Quadratic-linear schemes for the fractional integro-differential equations. The convergences of the presented numerical schemes are also established. The discussed schemes successfully validate the numerical results. The convergence of the numerical schemes is also discussed and validated through numerical results. It is observed that the Quadratic-linear scheme S3 performs comparatively better to schemes S1 and S2 for the test example 5.1. However, it appears that the scheme S2 performs comparatively good than schemes S1 and S3 for test example 5.2. Schemes S2 and S3 achieve second order convergence in both the test examples. In the third test example 5.3, the fractional integro-differential equation having exact solution as fractional power of $x$ is considered and numerical schemes are performed. And, in this case, it is observed that the scheme S1 achieves better convergence order than the schemes S2 and S3. Thus from these test examples, it is concluded that all the schemes perform well and provides accurate numerical results. It is also observed that the performance of the schemes depends on the nature of the problem. Here, the developed schemes are discussed only for the linear fractional integro-differential equations. However, similar study could be performed for the nonlinear fractional integro-differential equations in future.

## Acknowledgements

The authors sincerely thank to the reviewer for the useful comments and valuable suggestions toward the improvement of the manuscript. The first author acknowledges the financial supports from the university grant commission, New Delhi, India under the SRF schemes.